\documentclass[12pt]{amsart}
\usepackage{amsmath,amssymb,amscd,graphicx}
\usepackage{palatino}
\theoremstyle{plain}

\theoremstyle{definition}

\begin{document}
\title[Computing $L$-functions
with large conductor]{Computing $L$-functions\\
with large conductor}
\author{Jeffrey Stopple}
\email{stopple@math.ucsb.edu}
\maketitle

\begin{abstract}
An algorithm is given to efficiently compute
$L$-functions with large conductor in a restricted range of the critical strip.  Examples are included for about 21000 dihedral Galois representations
with conductor near $10^7$. The data shows good agreement with a symplectic
random matrix model.
\end{abstract}

\section{Introduction.}

In \cite{OS}, Odlyzko and Sch\"{o}nhage developed an algorithm to compute the Riemann zeta function $\zeta(s)$ efficiently for values of $s$ very high up in the critical strip.  Their method depends on precomputation of Taylor series expansions of $\zeta(s)$ at regularly spaced points, which in turn can be done efficiently by a clever application of the Fast Fourier Transform.  Rumely later implemented a version of this for Dirichlet $L$-functions in \cite{Rum}.

In analytic number theory it is often the case that there is a symmetry between what one can prove for large values of $t=\text{Im}(s)$  in the critical strip, and what one can prove for $L$-functions with large conductor $q$.  With this in mind, it seems reasonable to try to find an  algorithm to efficiently compute values of $L$-functions with very large conductor.

Roughly speaking, we want to view the $L$-function as a Mellin transform of an automorphic form $f$, and split the integral at the symmetry point.  This gives the extended $L$-function as an infinite sum 
\[
\Lambda(s,f)=\sum_n a(n)\{G(s,2\pi n/q^{1/2})+G(1-s,2\pi n/q^{1/2})\}
\]
of values of an incomplete Gamma function
\[
G(s,x)=x^{-s}\Gamma(s,x)=\int_1^\infty \exp(-yx)y^s \frac{dy}{y}.
\]
It was already observed in \cite {W} that the rapid decay of $G(s,x)$ as $x\to\infty$ means that we can truncate the infinite series at about $O(q^{1/2})$ terms.  On the other hand, the exponential decay as $t=\text{Im}(s)$ increases will cause a loss of precision if we want to compute $L$-function values for large $t$.  This can be fixed by moving the contour integral so that the $x$ parameter has complex values, as observed in \cite{MW} and \cite {LO} and implemented in \cite{Rub}.

The  algorithm is based on computing Taylor series expansions of the function $G(s,x)$ in the \emph{second} variable.  Differentiation under the integral sign and integration by parts gives a recursion which allows the computations of derivatives of $G(s,x)$ at very little cost.  The rapid decay of $G(s,x)$ as $x$ increases will imply that, instead of taking equally spaced points for the centers of the Taylor expansions, the sequence of points can grow exponentially.  This, in turn, will imply that we need very few Taylor expansions to compute efficiently, in fact $O(\log(q))$ expansions with $O(\log(q))$ terms in each.

However, too much change in $t$ will require that we increase the phase in the $x$ parameters, and then all Taylor expansions must be recomputed.  For that reason, we will only consider $|t|\le 1$ in this paper, although the method will generalize to any bounded region in the critical strip.  For further simplification we will consider only $s=1/2+it$ on the critical line, although this, too, is not crucial.

If we want to get the most benefit of the precomputations it makes sense to compute zeros of lots of $L$-functions with the same conductor simultaneously.  Here the power of the Fast Fourier Transform can again play a role, if all of the $L$-functions arise from the same abelian group structure.  In this paper we will focus primarily on $L$-functions attached to characters $\phi$ on the ideal class group of a complex quadratic field with discriminant $-q$.  The corresponding automorphic forms are weight one forms, which are cusp forms if the character $\phi$ is not a genus character.  The corresponding Galois representations are dihedral.  (The final section of the paper, however, includes an example of an elliptic curve $L$-function with large conductor, and of an $L$-function for a quadratic Dirichlet character for large discriminant.)

With this in mind we now explain the algorithm in a little more detail.  It will be convenient to use the correspondence between ideal classes and binary quadratic forms $Q$ of discriminant $-q$.  Corresponding to a character $\phi$ we have a theta function
\[
\Theta(z,\phi)=\sum_{[Q]}\phi(Q)\sum_{n>0}r_Q(n)\exp(2\pi inz)\overset{\text{def.}}=\sum_{n>0}r_\phi(n)\exp(2\pi inz).
\]
Then
\begin{multline*}
\Lambda(s,\phi)\overset{\text{def.}}=(q^{1/2}/2\pi)^s\Gamma(s)L(s,\phi)=\\
\sum_{n>0}r_\phi(n)\left\{G(s,2\pi n/q^{1/2})+G(1-s,2\pi n/q^{1/2})\right\}.
\end{multline*}
(See \cite[\S12.4]{Ibook} or \cite[\S22.3]{IK}.)
For purposes of locating zeros on the critical line it is better to compute the Hardy function $Z(t,\phi)$ defined for $t>0$ by
\[
\theta(t,\phi)= t \log(q^{1/2}/(2\pi)) + \arg(\Gamma(1/2 + i t))
\]
and, with $s=1/2+it$,
\begin{align*}
Z(t,\phi)=&\exp(i \theta(t,\phi))L(s,\phi)\\
=&\exp(i \theta(t,\phi))(2\pi/q^{1/2})^{s}\Gamma(s)^{-1}\times\\
&\sum_{n>0}r_\phi(n)\left\{G(s,2\pi n/q^{1/2})+G(1-s,2\pi n/q^{1/2})\right\}.
\end{align*}
We will truncate the infinite series after $N$ terms, and arrange the $n<N$ into intervals $I_j$, $j=1,2,\ldots T$ centered at points $x_j$ and with width $\Delta_j$.  For $n$ in the interval $I_j$ we compute a Taylor series expansion, truncated to $B$ terms
\[
G(s,2\pi n/q^{1/2})+G(1-s,2\pi n/q^{1/2})\approx
\sum_{k=0}^B  G_{j,k}(s)(2\pi n/q^{1/2}-x_j)^k,
\]
where
\begin{equation}\label{Eq:defgjk}
G_{j,k}(s)=\left\{G^{(k)}(s,x_j)+G^{(k)}(1-s,x_j)\right\}/k!.
\end{equation}
This gives  $Z(t,\phi)$ as
\begin{multline}\label{Eq:thm}
Z(t,\phi)\approx \exp(i \theta(t,\phi))(2\pi/q^{1/2})^{s}\Gamma(s)^{-1}\times\\
\sum_{j=1}^T\sum_{k=0}^B G_{j,k}(s)\Delta_j^k\sum_{n\in
I_j}r_\phi(n)\frac{(2\pi n/q^{1/2}-x_j)^k}{\Delta_j^k}.
\end{multline}
Here we have inserted canceling terms $\Delta_j^k$, $\Delta_j^{-k}$ to control the small size of $G_{j,k}(s)$ and the large size of $(2\pi n/q^{1/2}-x_j)^k$.

The point, then, is that the inner sum over $n\in I_j$ is independent of $s$, and so can be done as a precomputation.  Of course, this is done separately for each character $\phi$, but the coefficients $r_\phi(n)$ can be computed for all characters $\phi$ very efficiently via the Fast Fourier Transform once the representation numbers $r_Q(n)$ are known.

If we want to use this method for other kinds of $L$-functions, for example, the $L$-function of an elliptic curve, or a quadratic Dirichlet character, we forgo the advantage of using the Fast Fourier Transform.  However, we see below we can choose the parameters $x_j$ to be an integer times (a constant analogous to) $2\pi/q^{1/2}$, so that the precomputation can even be done with \emph{integer} arithmetic when the coefficients $r(n)$ are integers.  The $\Delta_j^k$ terms are simply omitted, and the (analog of the)  $(2\pi/q^{1/2})^k$ term is factored out of the inner sum.  The author would like to thank Kimberly Hopkins for help with this idea.

The remaining sections of the paper are as follows:
\begin{enumerate}
\item[2.] Review of properties of $G(s,x)$.
\item[3.] Truncation of the $L$-series after $N$ terms.
\item[4.] Arrangement of the Taylor expansions.  
\newline
The $T$ intervals $I_j=(x_j-\Delta_j/2,x_j+\Delta_j/2)$ are determined.
\item[5.] Truncation of the Taylor expansions after $B$ terms.
\item[6.] Implementation and examples.
\end{enumerate}

\section{Review of properties of $G(s,x)$.}

As above we define for $x>0$ and $s\in\mathbb C$
\[
G(s,x)=x^{-s}\Gamma(s,x)=\int_1^\infty \exp(-xy)y^s\frac{dy}{y}.
\]
Differentiating with respect to $x$ under the integral we see that
\begin{equation}\label{Eq:diffintegral}
\frac{d}{dx}G(s,x)=-\int_1^\infty\exp(-xy)y^{s+1}\frac{dy}{y}=-G(s+1,x).
\end{equation}
Integration by parts, on the other hand, gives
\begin{equation}\label{Eq:intbyparts}
G(s+1,x)=\frac{\exp(-x)}{x}+\frac{s}{x}G(s,x).
\end{equation}
Equations (\ref{Eq:diffintegral}) and (\ref{Eq:intbyparts}) give a nice recursive relation for all the derivatives $G^{(k)}(s,x)$ in terms of $G(s,x)$.

We will need bounds for $G(s,x)$ and its derivatives.  For $s$ in the critical strip, that is $0<\text{Re}(s)<1$, we have
\[
|G(s,x)| \le\int_1^\infty \exp(-xy)y^{\text{Re}(s)-1} dy.
\]
Change the variables by $u=y-1$ to get
\[
|G(s,x)|\le \exp(-x)\int_0^\infty \exp(-xu)(u+1)^{\text{Re}(s)-1} du.
\]
Now $u+1\ge 1$ and by hypothesis, $\text{Re}(s)-1<0$, so $(u+1)^{\text{Re}(s)-1}\le 1$, and therefore
\begin{equation}\label{Eq:biggbound}
|G(s,x)|\le \exp(-x)\int_0^\infty \exp(-xu)du=\frac{\exp(-x)}{x}.
\end{equation}
To estimate the derivatives $G^{(k)}(s,x)=(-1)^kG(s+k,x)$ this method does not apply, since we only have $\text{Re}(s)+k-1<k$.  Instead we can use the Cauchy formula for derivatives
\[
\frac{|G^{(k)}(s,x)|}{k!}\le \frac{M_R}{R^k},
\]
where $M_R$ is a bound for $G(s,w)$ with $|w-x|=R$.  An estimate similar to (\ref{Eq:biggbound}) shows that
$|G(s,w)|\le \exp(-\text{Re}(w))/\text{Re}(w)$, which a calculus argument shows is maximized at $w=x-R$, that is,
\begin{equation}\label{Eq:dergbound}
\frac{|G^{(k)}(s,x)|}{k!}\le \frac{\exp(R-x)}{(x-R)R^k}
\end{equation}
for any $0<R<x$.

\section{Truncation of the $L$-series.}

To compute to $D$ digits of accuracy, we need enough terms $N$ in the $L$-series so that
\begin{equation}\label{Eq:trunc1}
2\sum_{n>N}r_\phi(n)|\text{Re}(G(s,2\pi n/q^{1/2}))|\le10^{-D}.
\end{equation}
Since the $r_\phi(n)$ are the coefficients of a weight one cusp form, $r_\phi(n)\ll n^{1/2}$.  Using (\ref{Eq:biggbound}),  the left side of (\ref{Eq:trunc1}) is
\begin{align*}
\ll& \sum_{n>N} (q/n)^{1/2}\exp(-2\pi n/q^{1/2}) \\
\approx &q^{1/2}\int_N^\infty y^{-1/2}\exp(-2\pi y/q^{1/2})dy\\
<&(q/N)^{1/2}\int_N^\infty\exp(-2\pi y/q^{1/2})dy\\
\ll &q/N^{1/2}\exp(-2\pi N/q^{1/2}).
\end{align*}
For this to be less than $10^{-D}$ we need
\[
\log(q)+D\log(10) <2\pi N/q^{1/2}+\log(N)/2;
\]
it suffices that
\begin{equation}\label{Eq:bignchoice}
N=q^{1/2}\log(q\cdot 10^D)/2\pi .
\end{equation}

\section{Arrangement of the Taylor expansions.}

We seek to find intervals of radius $\Delta_j$ centered at $x_j$ and bounding circles of radius $R_j$ so that
\[
0<\Delta_j<R_j<x_j.
\]
The tail of the Taylor expansion
\begin{align}
&\sum_{k=B}^\infty  G_{j,k}(s)(2\pi n/q^{1/2}-x_j)^k\notag\\
\intertext{is bounded, via (\ref{Eq:dergbound})  by}
&\sum_{k=B}^\infty \frac{\exp(R_j-x_j)}{x_j-R_j}(\Delta_j/R_j)^k
=\frac{\exp(R_j-x_j)}{x_j-R_j}\frac{(\Delta_j/R_j)^B}{1-\Delta_j/R_j},\label{Eq:taylortail}
\end{align}
where $G_{j,k}(s)$ is as in (\ref{Eq:defgjk}).

One might wish to minimize the total number of terms in all the Taylor expansions, but this seems intractable.  Instead we will simply choose the truncation parameter $B$ to be independent of $j$.  To make this happen we require that
\begin{gather*}
\Delta_j/R_j=\quad \text{(constant)}\quad1/c <1,\\
x_j-R_j=\quad\text{(constant)}\quad K>0.
\end{gather*}
Furthermore, in order that the endpoints of the intervals meet, we require that
\[
x_j+\Delta_j=x_{j+1}-\Delta_{j+1}.
\]
This gives
\begin{gather*}
x_j=R_j+K=c\Delta_j+K\\
x_{j+1}=R_{j+1}+K=c\Delta_{j+1}+K,
\end{gather*}
so
\begin{gather*}
x_{j+1}-x_j=c(\Delta_{j+1}-\Delta_j)\\
x_{j+1}-x_j=\Delta_j+\Delta_{j+1}.
\end{gather*}
Thus
\begin{gather*}
\Delta_{j+1}+\Delta_j=c(\Delta_{j+1}-\Delta_j)\\
\Delta_{j+1}=\Delta_j(c+1)/(c-1).
\end{gather*}
Now
\[
x_{j+1}=x_j+\Delta_j+\Delta_{j+1}=x_j+\frac{2c}{c-1}\Delta_j
\]
and
\[
\Delta_j=R_j/c=(x_j-K)/c
\]
 implies that
\[
x_{j+1}=x_j(c+1)/(c-1)-2K/(c-1).
\]
For simplicity we choose the parameters $c=2$ and $K=x_1/2$ which gives
\begin{equation}\label{Eq:xjdef}
x_j=\frac{x_1}{2}\left(3^{j-1}+1\right),
\end{equation}
and
\begin{equation}\label{Eq:deltajdef}
\Delta_j=3^{j-1}\frac{x_1}{4}.
\end{equation}

The simplest choice for $x_1$ is $2\pi/q^{1/2}$.  This gives
\[
I_j=\frac{2\pi}{q^{1/2}}\cdot \left[\frac{3^{j-1}+2}{4},\frac{3^{j}+2}{4}\right].
\]
Observe that $I_1$ contains only the $n=1$ term, $I_2$ contains only the $n=2$ term, $I_3$ contains only $n=3,\ldots,7$, etc.
Thus for the very small values of  $j$  where the intervals contain fewer than $B$ terms one can do better computing each term instead of using the Taylor expansion.  This would improve the efficiency about 40\% for $q$ near $10^6$, 30\% near $10^9$,  and 20\%  near $10^{20}$.

We can now determine how many intervals $T$ are needed, since we need to compute terms in the series out to 
\[
n=N= q^{1/2}\log(q\cdot 10^D)/2\pi ,
\]
so
\[
\frac{2\pi N}{q^{1/2}} =\log(q\cdot 10^D).
\]
This requires that
\[
x_T+\Delta_T=\left(\frac{3^{T}}{2}+1\right)\frac{x_1}{2}>\log(q\cdot 10^D).
\]
With $x_1=2\pi/q^{1/2}$ we want
\[
2q^{1/2}\log(q\cdot 10^D)/\pi<3^T,
\]
or
\[
T>\left\{\log(2/\pi)+\log(q^{1/2}\cdot\log(q\cdot 10^D))\right\}/\log(3).
\]
With 
\[
\log(2/\pi)<0,\quad\text{ and }\quad
1/\log(3)<1
\]
 we can simply choose
\begin{equation}\label{Eq:bigtchoice}
T=\log(q^{1/2}\cdot\log(q\cdot 10^D))
\end{equation}

\section{Truncation of the Taylor expansions.}

Suppose we make error $\delta$ in computing
\[
G(s,2\pi n/q^{1/2})+G(1-s,2\pi n/q^{1/2}).
\]
Using the bound $n^{1/2}$ for $r_\phi(n)$ and the rough estimate $q^{1/2}$ for $N$,  the maximal error we make in the sum for $Z(t,\phi)$ is bounded by
\[
q^{-1/4}\sum_{n=1}^{q^{1/2}} n^{1/2}\cdot \delta \ll \delta \cdot q^{1/2}.
\]
Alternatively we can consider the standard error, assuming the errors in the terms are independent with standard deviation $\epsilon$.  Then the standard error in the sum is bounded by \cite{D}
\[
q^{-1/4}\left(\sum_{n=1}^{q^{1/2}} (n^{1/2}\cdot\epsilon)^2\right)^{1/2}\ll \epsilon \cdot q^{1/4}.
\]
We will assume the latter from now on.  We want $\epsilon\cdot q^{1/4} <10^{-D}$, or
\begin{equation}\label{Eq:epsdef}
 \epsilon<q^{-1/4}10^{-D},
\end{equation}
which will determine how many terms $B$ we need to take in each Taylor expansion.

It follows from the choice made in \S3 that
\[
\frac{\Delta_j}{R_j}=\frac{1}{2}\qquad\text{ and }\qquad x_j-R_j=x_1/2.
\]
Thus the tail (\ref{Eq:taylortail}) of the Taylor series reduces to
\begin{equation}\label{Eq:epsbnd}
\epsilon=\frac{2\exp(-x_1/2)}{x_1}2^{1-B}.
\end{equation}
Again, with $x_1=2\pi/q^{1/2}$ we combine (\ref{Eq:epsdef}) and (\ref{Eq:epsbnd}) to find that we require
\[
q^{1/2}\exp(-\pi/q^{1/2}) 2^{1-B}/\pi  < q^{-1/4}10^{-D}.
\]
Since 
\[
\exp(-\pi/q^{1/2})2/\pi <1
\]
we ignore it, and so we need
$
q^{3/4}10^D<2^B,
$
and let
\begin{equation}\label{Eq:bigbchoice}
B=1.5\cdot \log(q^{3/4}10^D)>\log_2(q^{3/4}10^D).
\end{equation}
Even under the assumption of maximal error we would still only need $B=O(\log(q10^D))$.

\section{Implementation and examples.}

To implement this algorithm requires the computation of reduced representatives of all the forms, as well as their coordinates in terms of the
generators of the cyclic factors of the class group.  This computation is clearly $O(q^{1/2})$.  

We then need to evaluate all the forms $ax^2+bxy+cy^2$  on a rectangular grid of integer lattice points $(x,y)$, large enough not to miss any representation of any integer $n<N$.  
Determining the dimensions of the grid is a Lagrange multipliers problem; we must maximize $g(x,y)=x$ (respectively $f(x,y)=y$) subject to the
constraint 
$ax^2+bxy+cy^2=N.$
One finds
\[
x\le2(Nc/q)^{1/2}\qquad\text{ and }\qquad y\le 2(Na/q)^{1/2}.
\]
Since the forms are reduced, the inequalities on $c$ and $a$ and our choice (\ref{Eq:bignchoice}) of $N$ imply
\[
x\le q^{1/4}\log(q\cdot 10^D)^{1/2}\qquad\text{ and }\qquad y \le (2/3^{1/4} )\log(q\cdot 10^D)^{1/2}.
\]
We have triple $(x^2,xy,y^2)$ for each lattice point, which are the rows of a $O(q^{1/4}\log(q\cdot 10^D))\times 3$ matrix.  Evaluation of the all the forms at all the lattice points consists of multiplying the above matrix by the $3\times O(q^{1/2})$ matrix of form data.  This is $O(q^{3/4}\log(q))$ multiplications (in terms of $q$, for a fixed number $D$ of digits of accuracy.)

The representation numbers $r_Q(n)$ are computed by brute force and ignorance; we look at each entry in the matrix product and increment the corresponding $r_Q(n)$, another $O(q^{3/4}\log(q))$ operations.

Fast Fourier Transform on a group of size $O(q^{1/2})$ is $O(q^{1/2}\log(q^{1/2}))$, and the function is vector valued with 
$N=O(q^{1/2}\log(q))$ entries, for a total computation of size $O(q\log(q)^2)$.

For each character $\phi$, the precomputation of the sums
\[
\sum_{n\in I_j}r_\phi(n)\frac{(2\pi n/q^{1/2}-x_j)^k}{\Delta_j^k}
\]
takes $N\cdot B$ operations, which is $O(q^{1/2}\log(q)^2)$ by (\ref{Eq:bignchoice}) and (\ref{Eq:bigbchoice}).  Subsequently, evaluations of 
$Z(t,\phi)$ using (\ref{Eq:thm}) require only
$T\cdot B$ operations, which is $O(\log(q)^2)$ by (\ref{Eq:bigtchoice}) and (\ref{Eq:bigbchoice}).

\begin{table}[!h]
\begin{center}
\begin{tabular}{ c c c c } 
-q&h(-q)&C(-q)& $\sharp$ $L$-functions\\ \hline
-10000003& 706& \{706\} &352\\
-10000004& 1648& \{412, 2, 2\}&820 \\
-10000007& 3660& \{3660\}&1829\\
-10000011& 816&\{204, 2, 2\}&404 \\
-10000015& 1134& \{1134\}&566\\
-10000019& 1275&  \{1275\}&637\\
-10000020& 928& \{232, 2, 2\}&460\\
-10000023&2064&\{258, 2, 2, 2\}&1024\\
-10000024&990&\{330, 3\}&494\\
-10000027&282&\{282\}&140\\
-10000031&5426&\{5426\}&NA\\
-10000036&876&\{876\}&437\\
-10000039&1912&\{956, 2\}&954\\
-10000043&618&\{618\}&308\\
-10000047&1512&\{756, 2\}&754\\
-10000051&742&\{742\}&370\\
-10000052&1692&\{846, 2\}&844\\
-10000055&3584&\{896, 2, 2\}&1788\\
-10000056&1480&\{370, 2, 2\}&736\\
-10000059&968&\{484, 2\}&482\\
-10000063&1722&\{1722\}&860\\
-10000072&724&\{724\}&361\\
-10000079&4147&\{4147\}&2073\\
-10000083&416&\{208, 2\}&206\\
-10000084&1364&\{682, 2\}&680\\
-10000087&1076&\{1076\}&537\\
-10000088&1512&\{126, 6, 2\}&752\\
-10000091&1382&\{1382\}&690\\
-10000095&2928&\{732, 2, 2\}&1460\\
-10000099&640&\{320, 2\}&318\\
&&&
\end{tabular}
\caption{Class group data for various discriminants}\label{Ta:data}
\end{center}
\end{table}

\begin{figure}[!t]
\includegraphics[scale=1, viewport=0 0 300 200,clip]{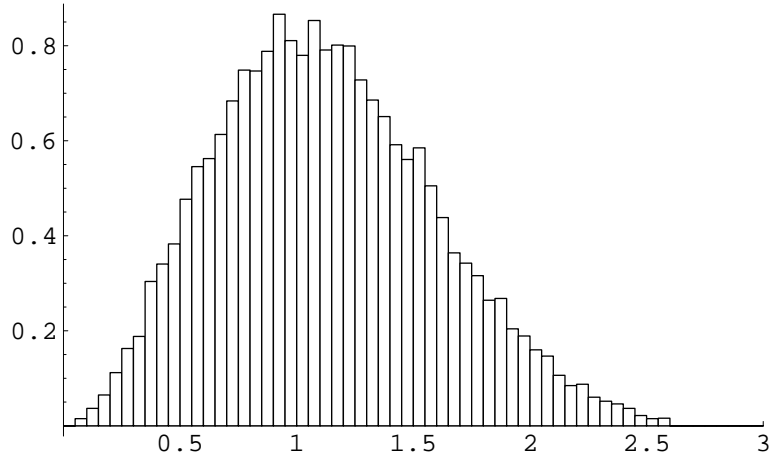}
\caption{Histogram data of first zero above $0.$}\label{F:first}
\end{figure}

This algorithm was implemented in \emph{Mathematica} and run on a 400 MHz. Apple Powerbook G4 under OS X 10.2\footnote{This paper was written in June 2003.}.  
\emph{Mathematica} is a good choice if one wants to re-invent the wheel as few times as possible.  
The incomplete Gamma
function $\Gamma(s,x)$ is available; it is computed via hypergeometric functions and continued fractions according to \cite[A.9.4]{MMA}.  
Also the Fast
Fourier Transform is supported in a sophisticated way \cite[A.9.4]{MMA}, via decomposition of the length of lists of data into prime factors.  
For large
factors, fast convolution methods are used.  
Since our data are real (the representation numbers $r_Q(n)$), the algorithm makes use of a real transform
method.  Non-cyclic abelian groups are handled automatically.      The
\emph{Mathematica} function \texttt{FindRoot} finds zeros of functions by a combination of damped
Newton's method, the secant method, and Brent's method \cite[A.9.4]{MMA}.

On the other hand, a specialized package for number theory, such as PARI \cite{PARI1} is preferable for computations involving binary quadratic
forms.   Fortunately, \emph{Mathematica}'s MathLink capability allows users to install their own C code, and this was done, essentially installing all the
PARI functions for computations with binary quadratic forms into \emph{Mathematica.}  (This was documented in \cite{PARI2}.)

The accuracy of the algorithm was checked by computing, to 15 digits, the zeros of some genus character $L$-functions.  By genus theory, these are products of Dirichlet $L$-functions attached to quadratic characters.  These zeros were first computed directly, writing the Dirichlet $L$-function as a linear combination of Hurwitz zeta functions (also supported in \emph{Mathematica}).  The zeros were also compared to the data in \cite{Rum}, and agreed to the number of digits given.

\begin{figure}
\includegraphics[scale=1, viewport=0 0 300 200,clip]{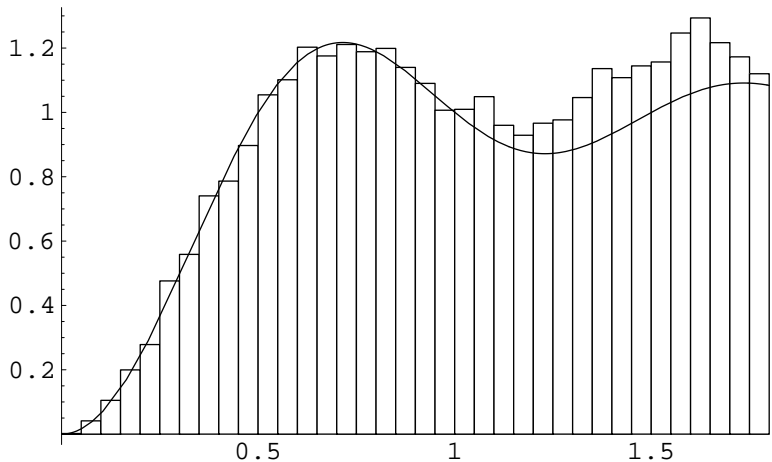}
\caption{$1$-level density vs. $1-\sin(2\pi x)/(2\pi x)$}\label{F:all}
\end{figure}

The algorithm was then run on 21336 $L$-functions with conductor near $q=10^7$.  
The Hardy function $Z(t,\phi)$ was evaluated, to six digits of accuracy, for $t$ between $0.$ and $1.$ in steps
of size $2\pi/(20\log(q))$, looking for sign changes.  When detected, \emph{Mathematica}'s \texttt{FindRoot} was used to find the zero.
This calculation took 137 hours, and found 35190 zeros below $t=1$.  
No attempt was made to prove GRH for these $L$-functions in this range, or to prove that all the zeros had been located.  The zeros themselves are available at \newline
\texttt{http://www.math.ucsb.edu/$\sim$stopple/}

\begin{figure}
\includegraphics[scale=1, viewport=0 0 300 300,clip]{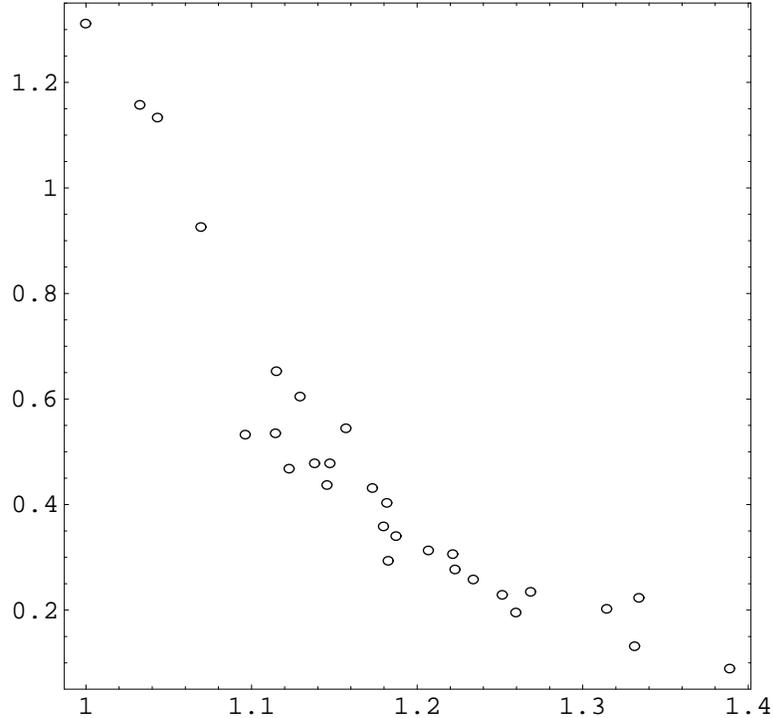}
\caption{Mean of first zero vs. $h(-q)/\sqrt{q}$.}\label{F:dotz}
\end{figure}

Table \ref{Ta:data} shows the discriminants, class numbers, class group structure, and number of $L$-functions (omitting real (genus) characters, and one of each complex conjugate pair.) Figure \ref{F:first} shows a histogram of of the lowest zero, renormalized by
$\tilde{\gamma}=\gamma\cdot \log(q)/(2\pi)$.  These seem to support a conjecture that the corresponding random matrix model is symplectic (see
\cite[Figure 5]{KS}).    The mean height above $0.$ is $1.13$, considerably larger than 
\[
0.78=\int_0^\infty t\nu_1(USp)(t)dt
\]
predicted by this model, due to the slow convergence as $q\to\infty$.  Figure \ref{F:all} compares the $1$-level density of $USp(\infty)$, that is, 
\[
1-\sin(2\pi x)/(2\pi x), 
\]
to the histogram data
\begin{gather*}
\frac{(\beta-\alpha)^{-1}}{\sharp \{\text{characters }\phi\}}\sum_\phi \sum_\gamma
\chi_{[\alpha,\beta)}(0.78/1.18\cdot\tilde{\gamma}),\\
\intertext{where}
\qquad [\alpha,\beta)=[0,.05),[.05,0.1),\dots, [1.75,1.80)
\end{gather*}   
The zeros $\tilde{\gamma}$ have been re-renormalized so that the mean above $0.$ of the first zero is $0.78$.

Finally, for each of the 29 discriminants considered, the mean above $0.$ of the first zero for the $L$-functions with that discriminant was compared to the relative size $h(-q)/\sqrt{q}$ of the class number.  
Figure \ref{F:dotz} shows the plot.  The discriminants with smaller class number tend to have a \emph{larger} mean first zero.  The
correlation between the two quantities is $-0.89$.   The connection between this phenomenon and the Chowla-Selberg formula is discussed in \S 3 of \cite{CFHR}.  

\begin{figure}
\includegraphics[scale=1, viewport=0 0 300 200,clip]{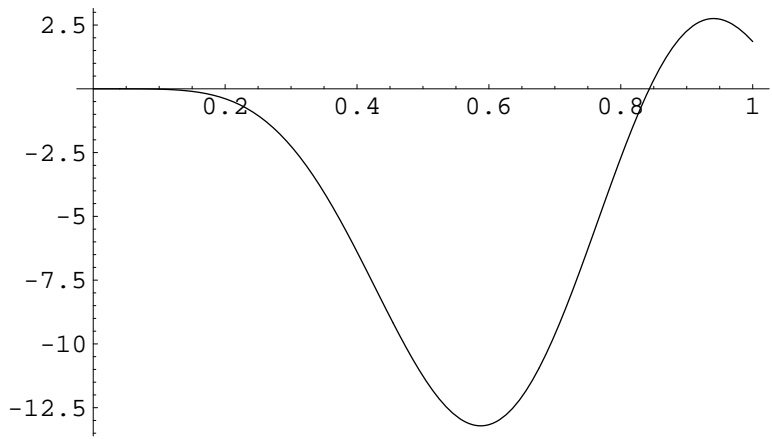}
\caption{Elliptic curve $L$-function with conjectural order $5$ zero at $s=1$ (i.e. $t=0$).}\label{F:plot763002}
\end{figure}

As mentioned in the introduction, this method can be adapted to compute other types of $L$-functions.  For an elliptic curve $L$-function, one gives up the advantage of the Fast Fourier Transform, but gains the advantage of doing the precomputations with integer arithmetic.  The elliptic curve
\[
E: x^3+y^3=763002
\]
has rank $5$ according to \cite{ER}, so by the Birch-Swinnerton Dyer conjectures the $L$ functions $L(s,E)$ should have a zero of order $5$ at $s=1$.  The conductor is $1746516156012$.  Figure \ref{F:plot763002} show a plot of the corresponding Hardy function $Z(t,E)$.

Low lying zeros of $L$-functions for quadratic Dirichlet characters $\chi_d$ are useful and scarce.  The algorithm described here is of no real use since the existence of such a zero can be ruled out by a single evaluation.  However, highly accurate values of zeros have been useful in the past, for example \cite{OR}.   With the parameter $D=110$, and discriminant $d=-175990483$ we can compute  $L(1/2+it,\chi_d)$ is zero,
\begin{multline*}
t=0.0004752439954201629008767557526752\\
   446841851348862432434240449427326484\\
   62812721184470556544512670480839630.
\end{multline*}
This was checked with the method of \cite{W}, and is accurate to more than $100$ digits.


\begin{thebibliography}{99}
\bibitem{CFHR} C. Bays, K. Ford, R. Hudson, and M. Rubinstein, \emph{Zeros of Dirichlet $L$-functions near the real axis and Chebyshev's bias}, J. Number Theory, \textbf{87}, 2001, pp. 54-76.
\bibitem{D} G. Dahlquist, Numerical Methods, Prentice Hall, 1974.
\bibitem{ER} N. Elkies and N. Rogers, \emph{Elliptic curves $x^3+y^3=k$ of high rank},in Algorithmic Number Theory, Lecture Notes in Comput. Sci., \textbf{3076}, Springer, Berlin, 2004, pp. 184-193.
\bibitem{FI} E. Fouvry and H. Iwaniec, \emph{Low-lying zeros of dihedral $L$-functions}, Duke J. Math., \textbf{116}, 2003, pp. 189-217.
\bibitem{Ibook} H. Iwaniec, Topics in Classical Automorphic Forms, Graduate Studies in Mathematics \textbf{17}, AMS, 1991.
\bibitem{IK} H. Iwaniec and E. Kowalski, Analytic Number Theory. AMS Colloquium Publications \textbf{53}, 2004.
\bibitem{KS} N. Katz, and P. Sarnak, \emph{Zeros of zeta functions and symmetry}, Bull. Amer. Math. Soc., \textbf{36}, no. 1, 1999, pp. 1-26.
\bibitem{LO} J. Lagarias, and A. Odlyzko, \emph{On computing Artin $L$-functions in the critical strip}, Math. Comp. \textbf{33}, 147, 1979,
pp.1081-1095.
\bibitem{MMA} S. Wolfram, The Mathematica Book, 4th ed., 1999.
\bibitem{MW} H. Montgomery, and P. Weinberger, \emph{Notes on small class numbers}, Acta Arith. \textbf{XXIV} 1974, pp. 529-542.
\bibitem{OR} A. Odlyzko and H. te Riele, \emph{Disproof of the Mertens' conjecture}, J. Reine Angew. Math. \textbf{357} (1985), pp. 138-160.
\bibitem{OS} A. Odlyzko and A. Sch\"{o}nhage,  \emph{Fast algorithms for multiple evaluations of the Riemann zeta function}, Trans. Amer. Math.
Soc. \textbf{309}, 2, 1988 pp. 797-809.
\bibitem{PARI1} PARI-GP. By C. Batut, D. Bernardi, H. Cohen , and M. Olivier, currently maintained by K. Belabas. \texttt{http://www.parigp-home.de}
\bibitem{PARI2} PARI-GP mailing list, available at\newline \texttt{http://www.parigp-home.de/lists/200212aav}
\bibitem{Rub} M. Rubinstein, \emph{Evidence for a spectral interpretation of the zeros of $L$-functions},  Ph.D. thesis, Princeton 1998.
\bibitem{Rum} R. Rumely, \emph{Numerical computations concerning the ERH}, Math. Comp., \textbf{61}, 203, 1993, pp. 415-440.
\bibitem{W} P. Weinberger, \emph{On small zeros of Dirichlet L-functions}, Math. Comp., \textbf{29} no. 129, 1975, pp. 319-328.

\end{thebibliography}
\end{document}